\input amstex.tex
\input amsppt.sty

\NoBlackBoxes

\define\Z{\Bbb Z}
\define\C{\Bbb C}
\define\R{\Bbb R}
\redefine\L{\Cal L}

\topmatter

\title
On a conjecture of Widom
\endtitle

\author Alexei Borodin and Alexei Novikov
\endauthor

\date August 25, 2005
\enddate

\abstract We prove a conjecture of H.~Widom stated in [W]
(math/0108008) about the reality of eigenvalues of certain
infinite matrices arising in asymptotic analysis of large Toeplitz
determinants. As a byproduct we obtain a new proof of
A.~Okounkov's formula for the (determinantal) correlation
functions of the Schur measures on partitions.
\endabstract

\endtopmatter

\head Introduction
\endhead

The main goal of this note is to prove that the spectrum of the operator $T$ in
$L^2(n,n+1,\dots)$, $n=0,1,\dots$, with matrix elements
$$
T_{pq}=\sum_{k\ge 1}\left(\frac {\phi_-}{\phi_+}\right)_{p+k}\left(\frac
{\phi_+}{\phi_-}\right)_{-q-k}
$$
where $(f)_k$ denotes the $k$th Fourier coefficient of the function $f(z)$ on
the unit circle $\{z\in\C:|z|=1\}$ and
$$
\phi_+(z)=e^{\gamma^+ z}\prod_{i\ge 1}\frac{
1+\beta_i^+z}{1-\alpha_i^+z},\qquad \phi_-(z)=e^{\gamma^-/z}\prod_{i\ge
1}\frac{ 1+\beta_i^-/z}{1-\alpha_i^-/z}
$$
for certain nonnegative parameters $\{\alpha_i^\pm\}$, $\{\beta_i^\pm\}$ and
$\gamma^\pm$ such that $\sum_i(\alpha_i^\pm+\beta_i^\pm)<\infty$ and
$\alpha_i^\pm,\beta_i^\pm\le \operatorname{const}<1$ for all $i\ge 1$, is real
and lies between 0 and 1.

This property of the operator $T$ turns out to be useful in the
asymptotic analysis of growing Toeplitz determinants with symbol
$\phi=\phi_+\phi_-$ or, more generally, of the distribution
functions of the largest parts of random partitions distributed
according to the associated Schur measure. The statement was
conjectured by H.~Widom in \cite{W} in the case when the
parameters $\alpha_i^+$, $\beta_i^-$, $\gamma^\pm$ vanish and the
number of nonzero $\alpha_i^-$'s and $\beta_i^+$'s is finite. Some
more special cases were proved in \cite{BDR}, \cite{W}. We refer
to \cite{W} for details and further references.

The key fact which allows us to prove the reality of spectrum is that the
resolvent $T(1-T)^{-1}$ of matrix $T$ after conjugation by a diagonal matrix
with plus-minus ones becomes totally positive. While proving this fact we
obtain as a byproduct a new proof of A.~Okounkov's formula for the correlation
functions of the Schur measures, see \cite{O} and also \cite{J}, \cite{R} for
other proofs. More specifically, we show that the Schur measure may be viewed
as an $L$-ensemble and explicitly compute the correlation kernel
$K=L(1+L)^{-1}$.

One of the authors (A.B.) was partially supported by the NSF grant DMS-0402047.

\head An operator identity
\endhead

In what follows $\Z_+$ denotes the set of nonnegative integers and all contour
integrals are taken over circles centered at the origin with radii close enough
to 1.

Let $\phi_+$ be a holomorphic function in the disc
$$
D_r=\{z\in\C: |z|\le r\}
$$
with radius $r>1$, and let $\phi_-$ be a holomorphic function outside the disc
$D_{1/r}$. In other words,
$$
\phi_{\pm}(z)=\sum_{n=0}^\infty {(\phi_{\pm})}_nz^{\pm n}, \qquad\quad
{(\phi_{\pm})}_n=O(r^{-n})\quad \text{as}\quad  n\to\infty.
$$
Assume also that $\phi_+$ and $\phi_-$ do not vanish on $D_r$ and
$\overline{\C}\setminus D_{1/r}$, respectively.

Consider an operator $L$ in $L^2(\Z_+)\oplus L^2(\Z_+)$ with matrix
$$
L=\bmatrix 0& A^t\\-B&0\endbmatrix,
$$
where the generating functions of matrix elements of $A$ and $B$ are
$$
\sum_{p,q\ge0} A_{pq}u^pv^q=\frac 1{u+v}\left(\frac{\phi_+(u)}{\phi_+(-v)}-1
\right),\quad \sum_{p,q\ge0} B_{pq}u^pv^q=\frac
1{u+v}\left(\frac{\phi_-(u^{-1})}{\phi_-(-v^{-1})}-1 \right).
$$
One can also write the matrix elements in terms of the contour integrals
$$
\gathered A_{pq}=\frac 1{(2\pi i)^2}\oint\oint\left(\frac
{\phi_+(u)}{\phi_+(-v)}-1\right)\frac{dudv}{(u+v)u^{p+1}v^{q+1}}\,,\\
B_{pq}=\frac 1{(2\pi i)^2}\oint\oint\left(\frac
{\phi_-(u^{-1})}{\phi_-(-v^{-1})}-1\right)\frac{dudv}{(u+v)u^{p+1}v^{q+1}}\,.
\endgathered
$$

These integral representations imply that $|A_{pq}|$ and $|B_{pq}|$ decay
faster than $\operatorname{const}\cdot x^{-p-q}$ for any $1<x<r$ as
$p+q\to\infty$. Thus, the sum of absolute values of matrix elements of $L$ is
finite, and $L$ is a trace class operator.

Note that the change $\phi_+(z)\leftrightarrow \phi_+^{-1}(-z)$ replaces $A$ by
$A^t$, the change $\phi_-(z)\leftrightarrow \phi_-^{-1}(-z)$ replaces $B$ by
$B^t$, and the change
$$
(\phi_+(u),\phi_-(v))\longleftrightarrow (\phi_-(u^{-1}),\phi_+(v^{-1}))
$$
switches $A$ and $B$. Thus, the switch $(A,B)\leftrightarrow (B^t,A^t)$ is
achieved by
$$
(\phi_+(u),\phi_-(v))\longleftrightarrow
(\phi_-^{-1}(-u^{-1}),\phi_+^{-1}(-v^{-1})).
$$

\proclaim{Theorem 1} Assume that the operator $1+A^tB$ is invertible
(equivalently, $1+BA^t$ is invertible). Then the operator $1+L$ is invertible,
and the matrix of the operator $K=L(1+L)^{-1}=1-(1+L)^{-1}$ has the form
$$
K=\bmatrix K_{11}&K_{12}\\K_{21}&K_{22}\endbmatrix=\bmatrix 1-(1+A^tB)^{-1}&
(1+A^tB)^{-1}A^t\\ -(1+BA^t)^{-1}B & 1-(1+BA^t)^{-1}
\endbmatrix
$$
where
$$
\gathered {(K_{11})}_{pq}=\delta_{pq}-\frac{(-1)^{p+q}}{(2\pi
i)^2}\oint\oint_{|zw|<1} \frac{\Phi(z,w)\,
dzdw}{(1-zw)\,z^{p+1}w^{q+1}}\,,\\
{(K_{12})}_{pq}=\frac{(-1)^{p}}{(2\pi i)^2}\oint\oint_{|zw|<1}
\frac{\Phi(z,w)\, w^q dzdw}{(1-zw)\,z^{p+1}}\,,\\
{(K_{21})}_{pq}=\frac{(-1)^{q+1}}{(2\pi i)^2}\oint\oint_{|zw|<1}
\frac{\Phi(z,w)\, w^{p}dzdw}{(1-zw)\,z^{q+1}}\,,\\
{(K_{22})}_{pq}=\frac{1}{(2\pi i)^2}\oint\oint_{|zw|<1} \frac{\Phi(z,w)\,
z^{p}w^q dzdw}{1-zw}\,,
\endgathered
$$
and
$$
\Phi(z,w)=\frac{\phi_-(z)\phi_+(w^{-1})}{\phi_+(z)\phi_-(w^{-1})}\,.
$$
\endproclaim
\demo{Proof} The proof of the first equality is straightforward. Thanks to the
symmetries mentioned before the statement of the theorem it suffices to prove
the integral formulas for $K_{11}$ and $K_{12}$. Let us start with $K_{11}$.

We need to show that $(1+A^tB)(1-K_{11})=1$ or, equivalently,
$A^tB(1-K_{11})=K_{11}$. Explicit computation gives
$$
\gather \bigl(A^tB(1-K_{11})\bigr)_{pq}=\sum_{l,m\ge0}
A_{lp}B_{lm}(1-K_{11})_{mq}\\ = \frac1{(2\pi i)^6}\oint\cdots\oint_{|zw|<1}
\sum_{l,m\ge0} \left(\frac{\phi_+(u_1)}{\phi_+(-v_1)}-1\right)\frac{du_1dv_1}
{(u_1+v_1)u_1^{l+1}v_1^{p+1}}\\
\times \left(\frac{\phi_-(u_2^{-1})}{\phi_-(-v_2^{-1})}-1\right)\frac{du_2dv_2}
{(u_2+v_2)u_2^{l+1}v_2^{m+1}}
\frac{\phi_-(z)\phi_+(w^{-1})}{\phi_+(z)\phi_-(w^{-1})}\,\frac{
dzdw}{(1-zw)\,(-z)^{m+1}(-w)^{q+1}}\\
=\frac1{(2\pi i)^6}\oint\cdots\oint_{|zw|<1,\,
|u_1u_2|>1,\,|v_2z|>1}\left(\frac{\phi_+(u_1)}{\phi_+(-v_1)}-1\right)
\frac{du_1dv_1}
{(u_1+v_1)v_1^{p+1}}\\
\times
\left(\frac{\phi_-(u_2^{-1})}{\phi_-(-v_2^{-1})}-1\right)\frac{du_2dv_2}
{(u_2+v_2)}
\, \frac{\phi_-(z)\phi_+(w^{-1})}{\phi_+(z)\phi_-(w^{-1})}\,\frac{
dzdw}{(-w)^{q+1}}\,\frac {(-1)}{(u_1u_2-1)(v_2z+1)(1-zw)}
\endgather
$$
where we imposed additional conditions $|u_1u_2|>1,\,|v_2z|>1$ on the
integration contours to ensure the convergence of two geometric series under
the integral.

We can immediately perform the integration over $u_2$ and $v_2$. Indeed, there
is only one simple pole $u_2=u_1^{-1}$ inside the $u_2$-contour and there is
only one simple pole $v_2=-z^{-1}$ inside the $v_2$-contour. Evaluating the
residues we obtain that the integral above equals
$$
\multline \frac1{(2\pi
i)^4}\oint\cdots\oint_{|zw|<1}\left(\frac{\phi_+(u_1)}{\phi_+(-v_1)}-1\right)
\frac{du_1dv_1}
{(u_1+v_1)v_1^{p+1}}\\
\times \left(\frac{\phi_-(u_1)}{\phi_-(z)}-1\right)\frac{1}{(z-u_1)} \,
\frac{\phi_-(z)\phi_+(w^{-1})}{\phi_+(z)\phi_-(w^{-1})}\,\frac{
dzdw}{(-w)^{q+1}}\,\frac {(-1)}{(1-zw)}\,.
\endmultline
$$
Let us choose the contours so that $|u_1|<|z|$ and open the parentheses
$\left(\frac{\phi_-(u_1)}{\phi_-(z)}-1\right)$ in the integral above. The
second term vanishes because it has no singularities inside the $u_1$-contour.
The first term has only one simple pole $z=u_1$ inside the $z$-contour, and the
integration over $z$ gives
$$
\multline \frac1{(2\pi
i)^3}\oint\oint\oint_{|u_1w|<1}\left(\frac{\phi_+(u_1)}{\phi_+(-v_1)}-1\right)
\frac{du_1dv_1}
{(u_1+v_1)v_1^{p+1}}\\
\times
\frac{\phi_-(u_1)\phi_+(w^{-1})}{\phi_+(u_1)\phi_-(w^{-1})}\,\frac{dw}
{(-w)^{q+1}}\,\frac
{(-1)}{(1-u_1w)}\,.
\endmultline
$$
Now let us choose the contours so that $|u_1|>|v_1|$ and open the parentheses
$\left(\frac{\phi_+(u_1)}{\phi_+(-v_1)}-1\right)$. The first term can be
integrated over $u_1$ --- there is only one pole $u_1=-v_1$ inside the
$u_1$-contour. Thus, the first term equals the corresponding residue, that is
$$
\frac{-1}{(2\pi i)^2}\oint\oint_{|v_1w|<1}
\frac{\phi_+(w^{-1})}{\phi_+(-v_1)}\,\frac{dv_1dw}
{(1+v_1w)v_1^{p+1}(-w)^{q+1}}\,.
$$
By deforming the $w$-contour to $\infty$ and picking up the residue at
$w=-v_1^{-1}$ we immediately see that this integral is equal to $\delta_{pq}$.

Now the second term is equal to
$$
\frac1{(2\pi i)^3}\oint\oint\oint_{|u_1w|<1,\,|u_1|>|v_1|} \frac{du_1dv_1}
{(u_1+v_1)v_1^{p+1}}\frac{\phi_-(u_1)\phi_+(w^{-1})}{\phi_+(u_1)\phi_-(w^{-1})}
\,\frac{dw} {(1-u_1w)(-w)^{q+1}}\,.
$$
Deforming the $v_1$-contour to $\infty$ we pick up the residue at $v_1=-u_1$
which gives
$$
\frac{-1}{(2\pi i)^2}\oint\oint_{|u_1w|<1}
\frac{\phi_-(u_1)\phi_+(w^{-1})}{\phi_+(u_1)\phi_-(w^{-1})} \,\frac{du_1dw}
{(1-u_1w)(-u_1)^{p+1}(-w)^{q+1}}\,,
$$
and this is exactly $(K_{11})_{pq}-\delta_{pq}$. The proof of the formula for
$K_{11}$ is complete.

In order to prove the formula for $K_{12}$ we need to show that
$(1+A^tB)K_{12}=A^t$ or $A^tBK_{12}=A^t-K_{12}$. The computation of
$A^tBK_{12}$ literally follows the above arguments for $K_{11}$ and leads to
the sum of two terms:
$$
\multline \frac{1}{(2\pi i)^2}\oint\oint_{|v_1w|<1}
\frac{\phi_+(w^{-1})}{\phi_+(-v_1)}\,\frac{w^q dv_1dw} {(1+v_1w)v_1^{p+1}}\\ +
\frac{1}{(2\pi i)^2}\oint\oint_{|u_1w|<1}
\frac{\phi_-(u_1)\phi_+(w^{-1})}{\phi_+(u_1)\phi_-(w^{-1})} \,\frac{w^q du_1dw}
{(1-u_1w)(-u_1)^{p+1}}\,.
\endmultline
$$

The first term is immediately seen to be equal to $A_{qp}$ and the second one
is exactly $-(K_{12})_{pq}$. \qed
\enddemo

\head Schur functions and Schur measures
\endhead

We refer the reader to \cite{M}, \cite{S} for general information on partitions
and symmetric functions.

Let $\Lambda$ be the algebra of symmetric functions. It can be viewed as the
algebra of polynomials in countably many indeterminates
$\Lambda=\C[h_1,h_2,\dots]$ where the indeterminates $h_k$ are the complete
homogeneous symmetric functions of degree $k$. We also agree that $h_0=1$ and
$h_{-k}=0$ for $k<0$.

The Schur symmetric functions $s_\lambda$ are parameterized by partitions
$\lambda$ and are expressed through $h_k$'s by the Jacobi-Trudi formula
$$
s_\lambda=\det[h_{\lambda_i-i+j}]_{i,j=1}^N
$$
where $N$ is any number greater or equal to the number of nonzero parts of
$\lambda$. The Schur functions form a linear basis in $\Lambda$.

An algebra homomorphism $\pi:\Lambda\to \C$ is uniquely determined by its
values on $h_k$'s, or by generating series of these values
$$
H^\pi(z)=\sum_{z=0}^\infty \pi(h_n)z^n.
$$

Recall that a sequence $\{a_n\}_{n=0}^\infty$ is called {\it totally positive}
if all minors of the matrix $[a_{i-j}]_{i,j\ge 0}$ are nonnegative. Here all
$a_{-k}$ for $k>0$ are assumed to be equal to zero. We will only consider
totally positive sequences with $a_0=1$; clearly, multiplication of all members
of a sequence by the same positive number does not affect total positivity.

The following statement was independently proved by
Aissen-Edrei-Schoenberg-Whitney in 1951 \cite{AESW}, \cite{E}, and
by Thoma in 1964 \cite{T}. An excellent exposition of deep
relations of this result to representation theory of the infinite
symmetric group can be found in Kerov's book \cite{K}.

\proclaim{Theorem 2} A sequence $\{a_n\}_{n=0}^\infty$, $a_0=1$, is totally
positive if and only if its generating series has the form
$$
\sum_{n=0}^\infty a_nz^n=e^{\gamma z}\frac{\prod_{i\ge 1}
(1+\beta_iz)}{\prod_{i\ge 1}(1-\alpha_iz)}=:F(\alpha,\beta,\gamma)
$$
for certain nonnegative parameters $\{\alpha_i\}$, $\{\beta_i\}$ and $\gamma$
such that $\sum_i(\alpha_i+\beta_i)<\infty$.

Equivalently, an algebra homomorphism $\pi:\Lambda\to \C$ takes
nonnegative values on all Schur functions if and only if the
sequence $\{\pi(h_n)\}_{n\ge 0}$ is totally positive, that is,
$H^\pi(z)=F(\alpha,\beta,\gamma)$ for a suitable choice of
parameters $(\alpha,\beta,\gamma)$.
\endproclaim

We will call a specialization $\pi:\Lambda\to\C$ {\it positive} if
$\pi(s_\lambda)\ge 0$ for all partitions $\lambda$. Thus, the theorem above may
be viewed as a classification of all positive specializations of the algebra of
symmetric functions.

Let us now consider specializations $\pi_+$ and $\pi_-$ of $\Lambda$ such that
$$
H^{\pi_+}(z)=\phi_+(z),\qquad H^{\pi_-}(z)=\phi_-(z^{-1})
$$
for the holomorphic functions $\phi^\pm$ of the previous section, and set
$s_{\lambda}^{\pm}:=\pi_{\pm}(s_\lambda)$.

Following A.~Okounkov \cite{O} assign to any partition $\lambda$ the following
(generally speaking, complex) weight
$$
P\{\lambda\}=\frac{s_\lambda^+s_\lambda^-}{Z}\,,\qquad Z=\exp{\sum_{k\ge 1} k
\bigl(\ln \phi_+(z)\bigr)_k \bigl(\ln \phi_-(z)\bigr)_k}.
$$
One can show that $\sum_{\lambda}P\{\lambda\}$ is an absolutely convergent
series with sum equal to 1. The distribution $P$ is called the {\it Schur
measure}.

Theorem 1 proved in the previous section yields a new proof of the
determinantal formula for the correlation functions of the Schur measure.

The following statement was proved in \cite{O}; other proofs can be found in
\cite{J}, \cite{R}.

For any partition $\lambda$ denote by $\L(\lambda)$ the infinite subset
$\bigl\{\lambda_i-i\bigr\}_{i=1}^\infty$ of $\Z$.

\proclaim{Theorem 3} For any $x_1,\dots,x_n\in\Z$
$$
\sum_{\lambda: \L(\lambda)\supset \{x_1,\dots,x_n\}}P\{\lambda\}= \det[\Cal
K(x_i,x_j)]_{i,j=1}^n
$$
where
$$
\Cal K(x,y)=\frac{1}{(2\pi
i)^2}\oint\oint_{|zw|<1}\frac{\phi_-(z)\phi_+(w^{-1})}{\phi_+(z)\phi_-(w^{-1})}\,
\frac{z^{p}w^q dzdw}{1-zw}\,.
$$
\endproclaim

\demo{Proof} We will use the material of Appendix in \cite{BOO}.
Let $\lambda$ be a partition and $(p_1,\dots,p_d\mid
q_1,\dots,q_d)$ be its Frobenius coordinates. Using the well-known
Giambelli formula $s_\lambda=\det[s_{(p_i|q_j)}]_{i,j=1}^d$, see
e.g. \cite{M, Ex.~I.3.9}, we obtain
$$
P\{\lambda\}=\frac{s_{\lambda}^+s_{\lambda}^-}{Z}=Z^{-1}
\det\left[s_{(p_i|q_j)}^+\right]_{i,j=1}^d\det\left[s_{(p_i|q_j)}^-
\right]_{i,j=1}^d.
$$
Comparing the following formula for the generating series of the hook Schur
functions
$$
1+(u+v)\sum_{p,q\ge 0} s_{(p|q)}u^pv^q=\frac {H(u)}{H(-v)},\qquad
H(z)=\sum_{n=0}^\infty h_nz^n,
$$
see e.g. \cite{M, Ex.~I.3.14}, with the definition of matrices $A$
and $B$ in the previous section we see that
$$
s_{(p|q)}^+=A_{pq},\qquad s_{(p|q)}^-=B_{pq}, \qquad p,q\ge 0.
$$
Thus,
$$
P\{\lambda\}=Z^{-1}\cdot \det [A_{p_iq_j}]_{i,j=1}^d \det
[B_{p_iq_j}]_{i,j=1}^d
$$
is, up to a constant, the value of the symmetric minor of the matrix $L$ from
the previous section with $2d$ rows and $2d$ columns marked by
$(q_1,\dots,q_d\mid p_1,\dots,p_d)$. This means that the Schur measure
interpreted through the Frobenius coordinates of partitions defines a
determinantal point process on $\Z_+\sqcup\,\Z_+$, and its correlation kernel
is given by $K=L(1+L)^{-1}$. (The operator $1+L$ is invertible because
$\det(1+L)=Z\ne 0$.)

Finally, to pass from Frobenius coordinates $(p_1,\dots,p_d\mid
q_1,\dots,q_d)\subset \Z_+\sqcup\,\Z_+$ to $\L(\lambda)\subset\Z$ we can use
the complementation principle, see \cite{\S A.3, BOO}, thanks to the Frobenius
lemma
$$
\L(\lambda)=\{p_1,\dots,p_d\}\sqcup\Bigl(\{-1,-2,-3,\dots\}\setminus
\{-q_1-1,\dots,-q_d-1\}\Bigr),
$$
see e.g. \cite{M, I.1(1.7)}. It is readily seen that the operation
$\triangle$ of \cite{\S A.3, BOO} transforms the kernel $K$ of
Theorem 1 to the kernel $\Cal K$ of Theorem 3 up to a conjugation
by a diagonal matrix of plus-minus 1's. Since such conjugation
does not change the determinants $\det[\Cal
K(x_i,x_j)]_{i,j=1}^n$, the proof is complete.\qed
\enddemo

\head Total positivity and Widom's conjecture
\endhead

In this section we assume that the functions $\phi_+(z)$ and $\phi_-(z^{-1})$
are generating functions of totally positive sequences, that is,
$$
\phi_+(z)=e^{\gamma^+ z}\frac{\prod_{i\ge 1} (1+\beta_i^+z)}{\prod_{i\ge
1}(1-\alpha_i^+z)},\qquad \phi_-(z)=e^{\gamma^-/z}\frac{\prod_{i\ge 1}
(1+\beta_i^-/z)}{\prod_{i\ge 1}(1-\alpha_i^-/z)}
$$
for certain nonnegative parameters $\{\alpha_i^\pm\}$, $\{\beta_i^\pm\}$ and
$\gamma^\pm$ such that $\sum_i(\alpha_i^\pm+\beta_i^\pm)<\infty$ and
$\alpha_i^\pm,\beta_i^\pm <r^{-1}<1$ for all $i\ge 1$.

Let $P_n$ be the projection operator in $L^2(\Z_+)$ which projects onto $n$
first basis vectors:
$$
P_n:(x_0,x_1,x_2,\dots)\mapsto (x_0,x_1,\dots,x_{n-1},0,0,\dots).
$$

\proclaim{Theorem 4} Under the above assumption, for any $n=1,2,\dots$ the
spectra of the operators $K^{(n)}_{11}=(1-P_n)K_{11}(1-P_n)$ and
$K^{(n)}_{22}=(1-P_n)K_{22}(1-P_n)$ are real and lie between 0 and 1.
\endproclaim

Let us start with a lemma.

\proclaim{Lemma 5} Let $T_n$ be a sequence of trace class operators in a
Hilbert space which converges to a (trace class) operator $T$ in trace norm.
Assume that $Sp(T_n)\subset S$ for a closed set $S\subset \C$ and all large
enough $n$. Then $Sp(T)\subset S\cup \{0\}$.
\endproclaim
\demo{Proof}  By virtue of \cite{Proposition A.11, BOO}, trace norm convergence
implies the convergence of entire functions in $z$
$$
\det(1+zT_n)\to \det(1+zT)
$$
uniformly on compact subsets of $\C$. Therefore, any zero of $\det(1+zT)$ is a
limit point of zeroes of $\det(1+zK_n)$. \qed
\enddemo
\demo{Proof of Theorem 4} Let us first show that $Sp(K_{11})$ and $Sp(K_{22})$
lie inside the interval $[0,1)$.

Since the specializations $\pi_\pm$ associated with $\phi_\pm$ are positive,
the matrices $A$ and $B$ are totally positive. Indeed, by Giambelli's formula
their minors are values of $\pi_\pm$ on appropriate Schur functions.
Furthermore, recall that $A_{pq}$ and $B_{pq}$ decay faster than
$\operatorname{const}\cdot x^{-p-q}$ for any $1<x<r$ as $p+q\to\infty$. This
means that the matrix $A^tB$ is also totally positive and its matrix elements
satisfy similar estimates. Consequently, $A^tB$ is a trace class operator, and
$\det(1+A^tB)-1$ is equal to the sum of all symmetric minors of $A^tB$, which
is nonnegative. Thus, the operator $1+A^tB$ is invertible.

The decay of the matrix elements of $A^tB$ shows that $P_nA^tBP_n$ converges to
$A^tB$ in trace norm as $n\to\infty$ (indeed, the trace norm of a matrix does
not exceed the sum of the absolute values of the matrix elements). It is well
known that the eigenvalues of a totally positive matrix are nonnegative, see
e.g. \cite{An, Corollary 6.6}. Hence, $Sp(P_nA^tBP_n)\subset \R_{\ge 0}$, and
by Lemma 5, $Sp(A^tB)\subset \R_{\ge 0}$. Therefore,
$Sp(K_{11})=A^tB(1+A^tB)^{-1}\subset [0,1)$. The argument for $K_{22}$ is very
similar.

To extend the argument to $K_{11}^{(n)}$, $K_{22}^{(n)}$ we need a linear
algebraic lemma.

\proclaim{Lemma 6} Let $C$ be a $(m+n)\times (m+n)$ matrix such that $1+C$ is
invertible, and let $D$ be the $n\times n$ lower right corner of $C(1+C)^{-1}$.
Then
$$
\det(1-D)=\sum_{X\subset \{1,\dots,m\}}
\frac{C\left(\matrix{X}\\{X}\endmatrix\right)}{\det(1+C)}
$$
(the sum above includes $X=\varnothing$ and
$C\binom{\varnothing}{\varnothing}=1$). Assume further that $\det(1-D)\ne 0$
and set $E=D(1-D)^{-1}$. Then any minor of $E$ is, up to a constant, a sum of
certain minors of the initial matrix $C$:
$$
\frac{E\left(\matrix{X}\\{Y}\endmatrix\right)}{\det(1+E)}=\sum_{Z\subset
\{1,\dots,m\}}\frac{C\left(\matrix {Z\sqcup X}\\{Z\sqcup
Y}\endmatrix\right)}{\det(1+C)}\,.
$$
\endproclaim
The proof is straightforward.

Now we can apply this statement to the totally positive matrices
$C=P_mA^tBP_m$, $m> n$. Then, clearly, $\det(1-D)>0$, and $E=D(1-D)^{-1}$ with
$$
D=(1-P_n)\frac{P_mA^tBP_m}{1+P_mA^tBP_m}(1-P_n),
$$
is totally positive. Since eigenvalues of totally positive matrices are
nonnegative, we obtain $Sp(D)\subset [0,1)$. As was mentioned earlier,
$P_mA^tBP_m\to A^tB$ in trace norm as $m\to\infty$. Hence, $D\to K_{11}^{(n)}$
in trace norm as $m\to\infty$, and Lemma 5 implies that
$Sp(K_{11}^{(n)})\subset [0,1]$. The case of $K_{22}^{(n)}$ is handled
similarly. \qed
\enddemo

Note now that by expanding $(1-zw)^{-1}=1+zw+(zw)^2+\dots$ in the integral
representation for $K_{11}$ we obtain
$$
\multline (K_{11})_{pq}=\delta_{pq}-(-1)^{p+q}\sum_{k\ge 0}\left(\frac
{\phi_-}{\phi_+}\right)_{p-k}\left(\frac {\phi_+}{\phi_-}\right)_{-q+k} \\ =
(-1)^{p+q}\sum_{k\ge 1}\left(\frac {\phi_-}{\phi_+}\right)_{p+k}\left(\frac
{\phi_+}{\phi_-}\right)_{-q-k}.
\endmultline
$$

H.~Widom in \cite{W} conjectured that the eigenvalues of the operator in
$L^2(n,n+1,\dots)$ with $(p,q)$-matrix element
$$
\sum_{k\ge 1}\left(\frac {\phi_-}{\phi_+}\right)_{p+k}\left(\frac
{\phi_+}{\phi_-}\right)_{-q-k}
$$
lie between 0 and 1 for the specific choice of
$$
\phi_+(z)=\prod_{i\ge 1}(1+r_iz),\qquad \phi_-(z)=\prod_{i\ge
1}(1-s_iz^{-1})^{-1}
$$
with nonnegative $r_i$ and $s_i$ only finitely many of which are nonzero;
$r_i,s_i< 1$. Clearly, this statement is an immediate corollary of Theorem 4.

\end